\documentclass[11pt]{article}
\usepackage{tikz}
\usetikzlibrary{calc,through,backgrounds}
\usepackage{multicol}

\usepackage{bm,amssymb,amsthm,amsfonts,amsmath,ytableau,latexsym,cite,color, psfrag,graphicx,ifpdf,pgfkeys,pgfopts,xcolor,extarrows,lipsum}

\newtheorem{theo}{Theorem}[section]

\newtheorem{lem} [theo]{Lemma}
\newtheorem{coro}[theo]{Corollary}
\newtheorem{prop}[theo]{Proposition}

\allowdisplaybreaks

\makeatletter \@addtoreset{equation}{section}
\@addtoreset{theo}{}\makeatother

\usepackage{hyperref}

\def\qed{\hfill \rule{4pt}{7pt}}
\def\pf{\noindent {\it Proof.} }

\def\inv{{\rm inv}}

\begin{document}

\begin{center}

{\bf INVERSION ARRANGEMENTS AND THE WEAK BRUHAT ORDER}\\[12pt]

\vskip 4mm

{\small NEIL J.Y. FAN}


\end{center}


\let\thefootnote\relax\footnotetext{This paper first appeared on the website www.paper.edu.cn in 2017.}

\noindent{A{\scriptsize BSTRACT}.}
For each permutation $w$,  we can construct a collection of hyperplanes $\mathcal{A}_w$ according to the inversions of $w$,   which is called the inversion hyperplane arrangement associated to $w$.  It was conjectured by Postnikov and confirmed by   Hultman,   Linusson,   Shareshian and   Sj\"{o}strand that the number of regions of $\mathcal{A}_w$ is less than or equal to the number of permutations below $w$ in the Bruhat order, with the equality holds if and only if $w$ avoids the four patterns 4231, 35142, 42513 and 351624. In this paper, we show that the number of regions of $\mathcal{A}_w$ is greater than or equal to the number of permutations below $w$  in the weak Bruhat order, with the equality holds if and only if $w$ avoids the patterns 231 and 312.




\section{Introduction}

Let $w=w_1w_2\cdots w_n$ be a permutation in $S_n$. The inversion hyperplane arrangement $\mathcal{A}_w$ associated to $w$ is the collection of hyperplanes $x_i-x_j=0$ in $\mathbb{R}^n$,  where $(i,j)$ is an inversion of $w$, i.e., $1\le i<j\le n$ and $w_i>w_j$. Denote by $re(w)$ the number of regions of $\mathcal{A}_w$, which are connected components of $\mathbb{R}^n\setminus\mathcal{A}_w $. In particular, if $w_0=n\, n-1\cdots 21$, then $\mathcal{A}_{w_0}$ is  the braid arrangement and $re(w_0)=n!$.

Let $br(w)$ (resp., $wk(w)$) denote the number of elements in the interval $[id,w]$ in the   Bruhat order (resp.,   weak Bruhat order), where $id=12\cdots n$ is the identity permutation.   Postnikov \cite{Postnikov} shows that $re(w)= br(w)$ when $w$ is a Grassmannian permutation and conjectures that
\begin{align}\label{rebr}
re(w)\le br(w),
\end{align}
 with the equality holds if and only if $w$ avoids the four patterns 4231, 35142, 42513 and 351624. In \cite{HLSS},
 Hultman,   Linusson,   Shareshian and   Sj\"{o}strand  confirm  this   conjecture of Postnikov. Furthermore, Hultman \cite{Hultman} provides a combinatorial criterion for $w\in W$ such that $re(w)= br(w)$, where $W$ is a finite Coxeter group.
 Oh, Postnikov and Yoo \cite{Oh} show  that the generating function for the regions of $\mathcal{A}_w$ coincides with the Poincar\'{e} polynomial of the interval $[id,w]$ in the Bruhat order if and only if $w$ avoids the patterns 3412 and 4231.

 The main result of this paper is

\begin{theo}\label{main}
 Let $w\in S_n$. Then we have
\begin{align}
re(w)\ge wk(w),
\end{align}
 with the equality holds if and only if $w$ avoids the patterns 231 and 312.
\end{theo}

Note that if a permutation avoids the patterns 231 and 312, it also avoids
the four patterns 4231, 35142, 42513 and 351624.  Combining the main result in \cite{HLSS} and Theorem \ref{main}, we obtain

\begin{coro}
For a permutation $w\in S_n$,  $wk(w)=br(w)$ if and only if $w$ avoids the patterns 231 and 312.
\end{coro}

Moreover, if $w$ avoids the patterns 231 and 312, then it avoids the patterns 3412 and 4231. Thus by \cite{Oh}, the generating function  for the regions of $\mathcal{A}_w$ coincides with the Poincar\'{e} polynomial of the interval $[id,w]$ in the (weak) Bruhat order.

To prove Theorem \ref{main}, we need to recall some connections of $re(w)$ with some other combinatorial structures.
For a permutation $w\in S_n$, the inversion graph $G_w$ associated to $w$ is a   graph with vertex set $\{1,\ldots,n\}$ and there is an edge between $i$ and $j$ if $(i,j)$ is an inversion of $w$. An acyclic orientation of   $G_w$ is an orientation of the edges of $G_w$ so that
the oriented graph has no directed cycles. Denoted by $ao(w)$ the number of acyclic orientations of $G_w$. It is well known that
\begin{align}\label{reao}
re(w)=ao(w),
\end{align}
see, for example,  Stanley \cite{Stanley}.

For a permutation $w\in S_n$, the south-west diagram $O_w$ of $w$ is a subset of $[n]\times [n]$ consisting of all pairs $(i,w_j)$ such that $i<j$ and $w_i<w_j$. A rook placement on a diagram (or board) $B$ is a set of cells (or ``rooks'') of $B$ such that no two rooks lie in the same row or column. For $B\subseteq[n]\times [n]$, denote $rk(B)$ by the number of placement of $n$ rooks on $([n]\times [n])\setminus B$, and let  $rk(w)=rk(O_w)$.
 Lewis and     Morales \cite{Lewis} show that
\begin{align}\label{aork}
ao(w)=rk(w).
\end{align}

\begin{figure}[!htb]
\setlength{\unitlength}{1.5mm}
\begin{center}
\begin{picture}(40,15)

\multiput(0,5)(3,0){5}{\circle*{0.4}}
\put(-.5,3){\footnotesize $1$}
\put(2.4,3){\footnotesize $2$}
\put(5.4,3){\footnotesize $3$}
\put(8.4,3){\footnotesize $4$}
\put(11.4,3){\footnotesize $5$}

\qbezier(0,5)(3,10.5)(6,5)\qbezier(3,5)(7.5,15)(12,5)
\qbezier(3,5)(6,10.5)(9,5)\qbezier(3,5)(4.5,7)(6,5)

\linethickness{0.05mm}
\multiput(30,3)(2,0){6}{\line (0,1){10}}
\multiput(30,3)(0,2){6}{\line (1,0){10}}

\linethickness{0.5mm}
\multiput(34,11)(0,2){2}{\line (1,0){6}}
\multiput(34,10.84)(2,0){4}{\line (0,1){2.32}}

\multiput(34,7)(0,2){2}{\line (1,0){4}}
\multiput(34,7)(0,2){2}{\line (1,0){4}}
\multiput(36,4.85)(2,0){2}{\line (0,1){4.32}}
\put(34,6.85){\line(0,1){2.32}}
\put(36,5){\line(1,0){2}}

\end{picture}
\end{center}
\vspace{-8mm}
\caption{The inversion graph $G_w$   and the south-west diagram $O_w$.}\label{fig}
\end{figure}
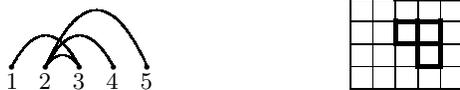

For example, let $w=25134$. Then $G_w$ and $O_w$ are displayed in Figure \ref{fig}. One can check that $ao(w)=rk(w)=16$.

For a permutation $w\in S_n$, the Lehmer code $c(w)$ of $w$ is the sequence $(c_1(w),c_2(w),\ldots,$ $c_n(w))$, where $c_i(w)=\#\{i<j\le n\,|\, w_j<w_i\}$.
Theorem \ref{main} follows from equations \eqref{reao}, \eqref{aork} and the following two propositions.
\begin{prop}\label{p1}
Let $w\in S_n$. Then we have
\begin{align}\label{eq1}
wk(w)\le \prod_{i=1}^n (c_i(w)+1),
\end{align}
with the equality holds if and only if $w$ avoids the pattern 231.
\end{prop}

\begin{prop}\label{p2}
Let $w\in S_n$. Then we have
\begin{align}\label{eq2}
\prod_{i=1}^n(c_i(w)+1)\le rk(w),
\end{align}
with the equality holds if and only if $w$ avoids the pattern 312.
\end{prop}

\section{The proof of Proposition \ref{p1}}

In this section, we give a proof of Proposition \ref{p1}. Let us first recall the Bruhat order and the (left) weak Bruhat order of the symmetric group.

The symmetric group $S_n$ can be viewed as a Coxeter group with generator set $S=\{s_1,\ldots,s_{n-1}\}$, where $s_i=(i,i+1)$ is an adjacent transposition.
A permutation $w=w_1\cdots w_n\in S_n$ can be viewed as a function on $[n]=\{1,\ldots,n\}$ such that $w_i=w(i)$.
Multiplications in $S_n$ are taken as composition of functions. Thus $s_iw$ is obtained by interchanging the elements $i$ and $i+1$ in $w$, whereas $ws_i$ is obtained by interchanging $w_i$ and $w_{i+1}$ in $w$. Denote
 \[I(w)=\{(i,j)\,|\,1\le i<j\le n, w_i>w_j\}\]
by the inversion set of $w$ and let $\inv(w)=|I(w)|$ be the number of inversions of $w$.

Let $u,v\in S_n$, we say that $u\le v$ in the Bruhat order if there exists a sequence of transpositions $(i_1,j_1),(i_2,j_2),\ldots,(i_k,j_k)$ such that
$v=u(i_1,j_1)(i_2,j_2)\cdots(i_k,j_k)$
and
 \[\inv(u(i_1,j_1)\cdots(i_r,j_r))=\inv(u)+r, \ \text{for}\   1\le r\le k.\]
  We say that $u\le v$ in the (left) weak Bruhat order if $v=s_1s_2\cdots s_k u$, where $s_i\in S$ and
$\inv(s_1s_2\cdots s_iu)=\inv(u)+i$, for $1\le i\le k$.
There is a well known criterion for comparing two permutations in the weak Bruhat order, see Bj\"{o}rner and   Brenti \cite{BB}.

\begin{lem}\label{wk}
Let $u,v\in S_n$.
Then $u\le v$ in the weak Bruhat order if and only if $I(u)\subseteq I(v)$.
\end{lem}

From Lemma \ref{wk} we can obtain the following  lemma directly.

\begin{lem}\label{wkc}
Let $u,v\in S_n$ such that $u\le v$ in the weak Bruhat order. Then $c(u)\le c(v)$ component-wisely, that is, $c_i(u)\le c_i(v)$ for $1\le i\le n$.
\end{lem}

Let $\pi=\pi_1\cdots\pi_k\in S_k$ be a permutation on $\{1,2,\ldots,k\}$. For $n\ge k$, we say that a permutation $w\in S_n$ contains the pattern $\pi$ if there exist $1\le i_1<\cdots<i_k\le n$ such that $w_{i_1}\cdots w_{i_k}$ is order isomorphic to $\pi_1\cdots\pi_k$. If $w$ does not contain the pattern $\pi$, we also say $w$ avoids the pattern $\pi$, or $w$ is $\pi$-avoiding.
Let
\[F(\Lambda_{w},q)=\sum_{u\le w}q^{\inv(u)},\]
where the sum is taken over $u\le w$ in the weak Bruhat order.
$F(\Lambda_{w},q)$ is called the Poincar\'{e} polynomial of the interval $[id,w]$ in the weak Bruhat order. The following theorem was also obtained by  Wei \cite{fan}. We provide a simple self-contained proof here.

\begin{theo}\label{wei}
Let $w\in S_n$ be a 231-avoiding permutation. Then
\[F(\Lambda_{w},q)=\prod_{i=1}^n [c_i(w)+1],\]
where $[k]=1+q+q^2+\cdots+q^{k-1}.$
\end{theo}

\pf Let $u\in S_n$ such that   $c(u)=(c_1(u),c_2(u),\ldots,c_n(u))$. According to Lemma \ref{wkc}, it suffices to show that if $c(u)< c(w)$ component-wisely then $u< w$ in the weak Bruhat order. Suppose to the contrary that $u\nless w$ in the weak Bruhat order. By Lemma \ref{wk}, we can choose an inversion $(i,j)\in I(u)\setminus I(w)$, namely, $u_i>u_j$ and $w_i<w_j$. Let $k\ (i<k\le j)$ be the smallest integer such that $w_i<w_k$. That is,
\begin{align*}
w&=w_1\cdots w_i\cdots w_k\cdots w_j\cdots w_n\\
u&=u_1\cdots u_i\cdots u_k\cdots u_j\cdots u_n.
\end{align*}
Note that since $w$ is 231-avoiding, all the elements $w_k,w_{k+1},\ldots,w_n$ are larger than $w_i$. And by the choice of $k$, all the elements $w_{i+1},w_{i+2},\ldots,w_{k-1}$ are smaller than $w_i$. Thus we obtain $c_i(w)=k-i-1$.

 We claim that all the elements $u_{i+1},u_{i+2},\ldots,u_{k-1}$ are smaller than $u_i$. If the claim holds, then we have $c_i(u)\ge k-i$, which is a contradiction.
To prove the claim, we first consider the element $u_{k-1}$. Since $w_{k-1}<w_i$ and $w_i<w_m$ for $m\ge k$, we see that $c_{k-1}(w)=0$. This implies  that $c_{k-1}(u)=0$. Thus $u_{k-1}<u_j<u_i$. Now consider the element $u_{k-2}$. Since $c_{k-2}(w)\le 1$, we see that $c_{k-2}(u)\le 1$. This leads to $u_{k-2}<u_i$. Since if $u_{k-2}>u_i$, then $c_{k-2}(u)\ge 2$. Continue this process, the claim follows. This completes the proof.
\qed

Now we are ready to give a proof of Proposition \ref{p1}.

\noindent{\it Proof of Proposition \ref{p1}}. According to Lemma \ref{wkc} and Theorem \ref{wei}, we need only to show that if $w$ contains a 231 pattern, then $wk(w)<\prod_{i=1}^n(c_i(w)+1).$ That is, we need to show that there exists $u\in S_n$ such that $c(u)< c(w)$ component-wisely but $u\nless w$ in the weak Bruhat order.

Assume that $w$ contains a 231 pattern. We claim that there exist $1\le i<j<j+1\le n$, such that $w_iw_jw_{j+1}$ forms a 231 pattern. In fact, choose a 231 pattern   in $w$ at random, say, $w_aw_bw_c$, where $a<b<c$ and $w_c<w_a<w_b$. If $c>b+1$, then consider the element $w_{c-1}$. There are  two cases. If $w_{c-1}>w_a$ then $w_aw_{c-1}w_c$ forms a claimed 231 pattern. If $w_{c-1}<w_a$, then we continue to consider the 231 pattern $w_aw_bw_{c-1}$ in the same manner. Since $c-b$ is finite, we see that the claim holds.
Choose the lexicographically smallest triple $(i,j,j+1)$ such that $w_iw_jw_{j+1}$ forms a 231 pattern, that is,
\[w=w_1\cdots w_i\cdots w_{j-1}w_jw_{j+1}w_{j+2}\cdots w_n.\]
Let
\[w'=w_1\cdots w_i\cdots w_{j-1}w_{j+1}w'',\]
where $w''$ is the sequence $w_{j}w_{j+2}\cdots w_n$ rearranged in the same relative order as  $w_{j+1}w_{j+2}$ $\cdots w_n$.
It is easy to check that $c(w')<c(w)$. In fact, only the $j$-th component of $c(w')$ is strictly smaller than that of $c(w)$, the other components are the same.  However, $I(w')\nsubseteq I(w)$, since $(i,j)\in I(w')\setminus I(w)$,  by Lemma \ref{wk}, $w'\nless w$ in the weak Bruhat order. This completes the proof.
\qed

For example, let $w=41382657$. Then $(1,4,5)$ is the lexicographically smallest triple such that $w_1w_4w_5$ forms a 231 pattern  and $c(w)=(3,0,1,4,0,1,0,0)$. Let $w'=4132w''=41325768$, where $w''=5768$ is obtained by rearrange $8657$ as the same relative order of $2657$. Then $c(w')=(3,0,1,0,0,1,0,0)<c(w)$ component-wisely but $w'\nless w$ in the weak Bruhat  order.

\section{The proof of Proposition \ref{p2}}

In this section, we provide a proof of Proposition \ref{p2}. Recall the definition of $O_w$ that, the number of cells $a_i(w)$  in the $i$-th row of $O_w$  is $a_i(w)=\#\{i<j\le n\,|\,w_j>w_i\}$. Clearly,  $a_i(w)=n-i-c_i(w)$.

\noindent{\it Proof of Proposition \ref{p2}}.
 Note that there is a rook in each row of  $\overline{O}_w=([n]\times [n])\setminus O_w$. To prove \eqref{eq2}, we claim that  the number of possible positions for  the rook in the $i$-th row of $\overline{O}_w$  is larger than or equal to $c_i(w)+1$. We make induction on the row number $i$. For $i=1$, since there are $n-c_1(w)-1$ cells in the first row of $O_w$, we see that there are $c_1(w)+1$ choices for the first rook, the claim holds.

Now assume that $i\ge 2$. For the   $i$-th row, there are $n-i-c_i(w)$ cells in $O_w$, so there are $c_i(w)+i$ cells in  $\overline{O}_w$. However, since each column contains exactly one rook, there are $i-1$ rooks in  $\overline{O}_w$ already, which may or may not occupy the columns that contain   cells in the $i$-th row of $O_w$. If all the first $i-1$ rooks do not occupy the columns that contain cells in the $i$-th row of $O_w$, then the possible positions for the $i$-th rook is $c_{i}(w)+1$. Otherwise, the possible positions for the $i$-th rook is strictly larger that $c_{i}(w)+1$.
 Therefore, the claim holds, and \eqref{eq2} follows.

We proceed to characterise the equality condition of \eqref{eq2}.  The crucial observation is that if $w$ is 312-avoiding, then $O_w$ is a right-justified Ferrers diagram, i.e., $a_1(w)\ge a_2(w)\ge\cdots\ge a_n(w)$ and the cells in the $i$-th row of $O_w$ are $(i,n-a_i(w)+1),\ldots,(i,n)$. To see this, suppose that there exists $a_i(w)<a_{i+1}(w)$  for some $1\le i<n-1$. Then we have $w_i>w_{i+1}$ and there must exist some $w_j$ such that $j>i+1$ and $w_{i+1}<w_j<w_i$. That is, $w_iw_{i+1}w_j$ forms a 312 pattern. Moreover, we claim that the cells in the $i$-th row of $O_w$ are $(i,n-a_i(w)+1),\ldots,(i,n)$. Suppose to the contrary  that in the $i$-th row, $(i,k)$ is a cell of $O_w$ but $(i,k+1)$ is not a cell of $O_w$. Then $k$ must appear to the right of $w_i$ and $w_i<k$, and so $k+1$ must appear to the left of $w_i$. Thus $(k+1)w_ik$ forms a 312 pattern, a contradiction.
According to the previous paragraph, we see that in the $i$-th row of $\overline{O}_w$, the number of possible positions for the $i$-th rook is exactly $c_i(w)+1$.

On the other hand, if $w$ contains a 312 pattern, then by applying similar arguments to the proof of Proposition \ref{p1}, we can choose $1\le i<i+1<j\le n$ such that $w_iw_{i+1}w_j$ forms a 312 pattern. Thus   $(i+1,w_j)$ is a cell of $\overline{O}_w$, and    $(i,w_j)$ is  not a cell  of $\overline{O}_w$. Therefore, the possible positions of the $(i+1)$-th rook in $\overline{O}_w$ is strictly larger than $c_{i+1}(w)+1$. This completes the proof. \qed

To conclude, we remark that in \cite{HLSS}, Hultman, Linusson,   Shareshian, and Sj\"{o}strand  construct a map $\phi$ to show the inequality \eqref{rebr}. In general, the map $\phi$ is injective but not surjective, so they ask a question to determine the image of $\phi$. Hultman \cite{Hultman} characterizes the condition under which the map $\phi$ is surjective for finite Coxeter groups.
We conjecture that the elements of the interval $[id,w]$ in the weak Bruhat order are   images of the map $\phi$. For the symmetric group, this conjecture was verified for $n\le 8$ by computer.

It is worth mentioning that the four patterns 4231, 35142, 42513 and 351624 have a geometric connection with Schubert varieties. More precisely, Gasharov and Reiner \cite{Gasharov} show that these four patterns can characterize   Schubert varieties in partial flag manifolds defined by inclusion.
On the other hand,   the Schubert variety $X_w$ indexed by $w$ is smooth if and only if $w$ avoids the patterns 3412 and 4231.
Therefore, it is natural to ask whether there is any connection between   geometry and the patterns 231 and 312.

\noindent{\bf Acknowledgements.}
This work was
supported by  the National Natural  Science Foundation of China.

\footnotesize{

N{\scriptsize EIL} J.Y. F{\scriptsize AN}, D{\scriptsize EPARTMENT OF} M{\scriptsize ATHEMATICS}, S{\scriptsize ICHUAN} U{\scriptsize NIVERSITY}, C{\scriptsize HENGDU} 610064, P.R. C{\scriptsize HINA.} Email address: fan@scu.edu.cn

}

\end{document}